\documentclass[a4paper]{article}
\usepackage{latexsym}
\usepackage{amsmath,amssymb,amsfonts}

\newcommand{\noi}{\noindent}

\newtheorem{prethm}{{\bf Theorem}}

\newtheorem{prepro}[prethm]{Proposition}

\newtheorem{prelem}[prethm]{Lemma}

\newtheorem{precor}[prethm]{Corollary}

\newtheorem{preremark}{{\bf Remark}}

\newtheorem{preexample}{{\bf Example}}

\newtheorem{preproof}{{\bf Proof.}}

\date{}
\title{\bf\ A conjecture on square roots of Laplacian and signless Laplacian eigenvalues of graphs }
\author{S. Akbari$^{\,\rm a,b}$,~~  E. Ghorbani$^{\,\rm a,b}$,~~  M. R. Oboudi$^{\,\rm a}$\\
{\small {\em $^{\rm a}$Department of Mathematical Sciences, Sharif University of Technology,  Tehran, Iran}}\\
{\small {\em $^{\rm b}$School of Mathematics, Institute for Research in Fundamental Sciences (IPM),}}\\
 {\small {\em P.O. Box 19395-5746, Tehran, Iran}}}
\begin{document}
\maketitle

Let $G$ be a graph with $n$ vertices, $m$ edges, the incidence
matrix $X$, and the directed incidence matrix $D$ for a given
orientation.
 $X$ and $D$ are $n$-by-$m$ matrices, $X$ is a $(0,
1)$-matrix, and $D$ a $(-1, 0, 1)$-matrix. Moreover
$$XX^\top=\Delta+A,~\hbox{and}~~DD^\top=\Delta-A,$$ where $A$ is the
adjacency matrix and $\Delta$ is the diagonal matrix whose diagonal
entries are the vertex degrees of $G$. Note that $\Delta-A$ is the
Laplacian matrix and $\Delta+A$ is the signless Laplacian matrix of
$G$. Thus the singular values of $X$ and $D$ are the square roots of
the signless Laplacian and the Laplacian eigenvalues  of $G$,
respectively. We denote the sum of all singular values of a matrix
$M$ by ${\cal E}(M)$. We have the following conjecture:

\noi{\bf Conjecture.} For any graph, ${\cal E}(D) \le {\cal E}(X).$

It is not hard to see that for bipartite graphs, ${\cal E}(D) ={\cal
E}(X).$ A computer aid search showed that the conjecture is true for
all graphs up to 9 vertices.  For  $k$-regular graphs with adjacency
eigenvalues $\lambda_1,\ldots,\lambda_n$ the conjecture is
equivalent to the following inequality:
$$\sqrt{k-\lambda_1}+\cdots+\sqrt{k-\lambda_n}\le \sqrt{k+\lambda_1}+\cdots+\sqrt{k+\lambda_n}.$$


\end{document}